\documentclass{article}
\usepackage{fullpage}
\def\C{{\rm C \kern-.48em\vrule width.06em height.57em depth-.02em \kern.48em}}
\def\Z{{{\rm Z}\kern-.28em{\rm Z}}}
\def\N{{{\rm I}\kern-.16em{\rm N}}}
\def\R{{{\rm I}\kern-.16em{\rm R}}}
\def\Re{\mathop{\rm Re}\nolimits}
\def\set#1{\langle #1 \rangle}
\def\belowrightarrow#1{{{{}\over\ #1\ }\kern-1.1em\to}}
\def\mystrut{\vphantom{\vrule width .1em height 1.em depth .7em}}
\def\goback#1{\setbox0\hbox{#1}\kern-\wd0 \relax}
\def\dbl#1{#1,#1}
\def\clm{\colon}
\def\y{(1-\lambda)}

\begin{document}

\title{Not all GKK $\tau$-matrices are stable}
\author{Olga Holtz \\  \small
Department of Mathematics \\
\small University of Wisconsin \\ \small Madison, Wisconsin 53706 U.S.A. \\
\small holtz@math.wisc.edu \normalsize}
\date{}
\maketitle

\begin{abstract} 
\noindent Hermitian positive definite, totally positive, and
nonsingular $M$-matrices enjoy many common properties, in particular
\begin{description}
\item{(A)} positivity of all principal minors,
\item{(B)} weak sign symmetry,
\item{(C)} eigenvalue monotonicity,
\item{(D)} positive stability.
\end{description}
The class of GKK matrices is defined by properties~(A) and~(B),
whereas the class of nonsingular $\tau$-matrices by~(A) and (C). 
It was conjectured that

\begin{description}
\item (A), (B) $\Longrightarrow$ (D) \hspace{0.2cm}  [D. Carlson, J. Res. Nat. Bur. Standards Sect. B 78 (1974) 1-2], 
\item (A), (C) $\Longrightarrow$ (D) \hspace{0.2cm} [G.M. Engel and H. Schneider, Linear and Multilinear Algebra 4 (1976) 155-176], 
\item (A), (B) $\Longrightarrow$ a property stronger than (D) \hspace{0.2cm} [R. Varga, Numerical Methods in Linear Algebra, 1978, pp.5-15], 
\item (A), (B), (C) $\Longrightarrow$ (D) \hspace{0.2cm} [Hershkowitz, Linear Algebra Appl. 171 (1992) 161-186].   
\end{description}

\noindent We describe a class of unstable GKK $\tau$-matrices, thus
disproving all four conjectures.
\end{abstract} 

\section{Definitions and notation}

Given a matrix $A\in \C^{n\times n}$, let $A(\alpha,\beta)$ denote the
submatrix of $A$ whose rows are indexed by $\alpha$ and columns by
$\beta$ ($\alpha$, $\beta\in \set{n}\clm =\{1,\ldots,n\}$) and let
$A[\alpha,\beta]$ denote $\det A(\alpha,\beta)$ if $\# \alpha=\#
\beta$ (where $\#$ stands for the cardinality of a set) with the
convention $A[\dbl\emptyset]\clm =1$. 

A matrix $A$ is called a {\em $P$-matrix \/} if $A[\alpha,\alpha]>0$
$\;\forall \alpha\subseteq \set{n}$. $A$ is {\em weakly
sign-symmetric\/} if
$$ A[\alpha,\beta] A[\beta,\alpha]\geq 0 \qquad \forall \alpha,\beta\in
\set{n}, \quad \# \alpha=\# \beta=\# \alpha \cup \beta -1.$$
Weakly sign-symmetric $P$-matrices are also called {$\!$\em GKK\/} 
after Gantmacher, Krein, and Kotelyansky.  It was proven by Gantmacher, Krein~\cite{1}, and
Carlson~\cite{2} that a $P$-matrix is GKK iff it satisfies the generalized
Hadamard-Fisher inequality
\begin{equation}
A[\dbl\alpha] A[\dbl\beta] \geq A[\dbl{\alpha\cup \beta}]
A[\dbl{\alpha \cap \beta}] \qquad \forall \alpha,\beta \subseteq
\set{n}. \label{HF}
\end{equation}
Carlson~\cite{3} conjectured that the GKK matrices are positive
stable, i.e., $\Re \lambda>0$ $\; \forall \lambda\in \sigma(A)$ (here
$\sigma(A)$ denotes, as usual, the spectrum of $A$), and showed his
conjecture to be true for $n\leq 4$.

Let $$l(A)\clm = \cases{ \min\{ \lambda \in \sigma(A) \cap \R \} & if
$\sigma(A) \cap \R\neq \emptyset$ \cr \infty & otherwise. }$$ $A$ is
called an $\omega$-matrix if it has eigenvalue monotonicity
$$ l(A(\alpha,\alpha))\leq l(A(\beta,\beta))<\infty \qquad \mbox{whenever}
\quad \emptyset \neq \beta \subseteq \alpha \subseteq \set{n}.$$ $A$ is a
$\tau$-matrix if, in addition, $l(A) \geq 0$. 

Engel and Schneider~\cite{ES} asked if nonsingular $\tau$-matrices or,
equivalently, $\omega$-matrices all whose principal minors are positive 
 (see Remark~3.7 in~\cite{ES}), are positive stable. Varga~\cite{V} 
conjectured even more than stability, viz.
$$ |\arg(\lambda - l(A))|\leq \frac{\pi}{2}-\frac{\pi}{n} \qquad
\forall \lambda \in \sigma(A). $$ This inequality was proven for
$n\leq 3$ by Varga (unpublished) and Hershkowitz and Berman~\cite{HB} and
for $n=4$ by Mehrmann~\cite{M}.  In his survey paper~\cite{H},
Hershkowitz posed the weaker conjecture that $\tau$-matrices that are
also GKK are stable.

Below we describe a class of GKK $\tau$-matrices which are not even
nonnegative stable, i.e., have eigenvalues with negative real part.
We construct Toeplitz Hessenberg matrices $A_{n,k,t}$ of size $n$ for
$k\in \N$ and $t\in \R$.  We show that $A_{n,k,t}$ is GKK for any
$t\in (0,1)$, a $\tau$-matrix if $n\leq 2k+2$ and $t\in(0,1)$ is
sufficiently small, and that $A_{2k+2,k,t}$ is unstable for sufficiently 
large $k$ and sufficiently small positive $t$. This provides a 
counterexample to the Hershkowitz conjecture and, therefore, to the 
Carlson, Engel and Schneider, and Varga conjectures as well.

In what follows, we shall use the following notation 
$$ p\clm q\clm=\cases{ \{p,p+1,\ldots,q\} & if $p\leq q$ \cr \emptyset &
                 otherwise} \qquad \forall p,q\in \N, \qquad \qquad
                 \quad x_+\clm=\cases{ x & if $x>0$ \cr 0 &
                 otherwise \cr} \qquad \forall x\in \R.$$

\section{Counterexample}

Given $k$, $n\in \N$, and $t\in (0,1)$, let $A_{n,k,t}$ be the following
Toeplitz Hessenberg matrix. If $n\leq k+1$, set
$$ A_{n,k,t}\clm=\left( \begin{array}{ccccc} 1 & 0 & \cdots & 0 & 0 \\
1 & 1 & \cdots & 0 & 0 \\ \vdots & \vdots & \ddots & \vdots & \vdots
\\ 0 & 0 & \cdots & 1 & 0 \\ 0 & 0 & \cdots & 1 & 1
\end{array}\right)_{n\times n}\goback{$n\times$}. $$ Otherwise let $$
A_{n,k,t}\clm=\left( \begin{array}{ccccccc} 1 & \overbrace{0 \; \cdots
\; 0 \; 0}^k & a^{k,t}_1 & a^{k,t}_2 & \cdots & a^{k,t}_{n-k-2} &
a^{k,t}_{n-k-1} \\ 1 & 1 \; \cdots \; 0 \; 0 & 0 & a^{k,t}_1 & \cdots
& a^{k,t}_{n-k-3} & a^{k,t}_{n-k-2} \\ \vdots & \! \vdots \ \; \ddots
\; \vdots \, \; \vdots & \vdots & \vdots & \ddots & \vdots & \vdots \\
0 & 0 \; \cdots \; 1 \; 0 & 0 & 0 & \cdots & a^{k,t}_1 & a^{k,t}_2 \\
0 & 0 \; \cdots \; 1 \; 1 & 0 & 0 & \cdots & 0 & a^{k,t}_1 \\ 0 & 0 \;
\cdots \; 0 \; 1 & 1 & 0 & \cdots & 0 & 0 \\ 0 & 0 \; \cdots \; 0 \; 0
& 1 & 1 & \cdots & 0 & 0 \\ \vdots & \! \vdots \ \; \ddots \; \vdots
\, \; \vdots & \vdots & \vdots & \ddots & \vdots & \vdots \\ 0 & 0 \;
\cdots \; 0 \; 0 & 0 & 0 & \cdots & 1 & 0 \\ 0 & 0 \; \cdots \; 0 \; 0
& 0 & 0 & \cdots & 1 & 1 \end{array} \right)_{n\times n} $$ where
$a^{k,t}_j$'s are chosen so that $A_{n,k,t}[\dbl{\set{k+j+1}}]=t^j$.
This definition makes sense for all $j=1,\ldots,n-k-1$. Indeed, the
expansion of $A_{n,k,t}[\dbl{\set{k+j+1}}]$ by the first row gives
\begin{eqnarray} && A_{n,k,t}[\dbl{\set{k+j+1}}](=\det A_{k+j+1,k,t}) = 
\nonumber \\ && \qquad A_{n,k,t}[\dbl{2\clm k+j+1}]+
\sum_{l=1}^j(-1)^{k+l} a^{k,t}_l A_{n,k,t} [\dbl{k+l+2\clm
k+j+1}]=\nonumber \\ && \qquad \det A_{k+j,k,t}+\sum_{l=1}^j(-1)^{k+l}
a^{k,t}_l \det A_{j-l,k,t} \label{!} \end{eqnarray} (recall that
$A_{n,k,t}[\dbl{\emptyset}]=1$, so the last term in the sum is well
defined).  As the coefficient of $a^{k,t}_j$ in~Eq.~(\ref{!}) is equal to
$(-1)^{k+j}$, the equation $A_{n,k,t}=s$ (linear in $a^{k,t}_j$) has a
solution for any right hand side $s$, in particular, for
$s\clm=t^j$. Since $A_{n,k,t}$ is Toeplitz, this implies
$A_{n,k,t}[\dbl{i\clm i+j-1}]=t^{(j-k-1)_+}$.

Show that the matrices $A_{n,k,t}$ are GKK for any $t\in (0,1)$. 
Since $A_{n,k,t}$ is Hessenberg, the submatrix
$A_{n,k,t}(\dbl{\set{n} \setminus i\clm i+j-1})$ is block upper
triangular if $1<i\leq i+j-1<n$, so
\begin{equation}
A_{n,k,t}[\dbl{\alpha\cup \beta}]=A_{n,k,t}[\dbl \alpha] A_{n,k,t}
[\dbl \beta] \qquad {\rm whenever} \qquad i<j-1 \;\; \hbox{\rm for all
} i\in \alpha, \; j \in \beta. \label{*}
\end{equation}
This shows that $A_{n,k,t}$ is a $P$-matrix. Moreover, since $0<t<1$
and $$ (x+y-k-1)_+ + (x+z-k-1)_+ \leq (x-k-1)_+ + (x+y+z-k-1)_+ \qquad
\forall x,y,z\geq 0, $$ we have
\begin{equation}
\begin{array}{l}
\mystrut A_{n,k,t}[\dbl{i\clm i+j-1}] \cdot A_{n,k,t}[\dbl{l\clm
l+m-1}]=\\ \mystrut \quad t^{(j-k-1)_+ + (m-k-1)_+}\geq
t^{(l+m-i-k-1)_+ + (i+j-l-k-1)_+}=\\ \mystrut \quad
A_{n,k,t}[\dbl{i\clm l+m-1}]\cdot A_{n,k,t}[\dbl{l\clm i+j-1}]
\end{array} \qquad {\rm if} \; l\leq i+j-1.
\label{**} \end{equation}
Together with~Eq.~(\ref{*}), Eq.~(\ref{**}) shows that $A_{n,k,t}$
satisfies~Eq.~(\ref{HF}) if $\alpha$, $\beta$ are sets of consecutive
integers.

To prove~Eq.~(\ref{HF}) in general, first make a definition. Call the
subsets $\alpha$, $\beta \subseteq \set{n}$ {\em separated\/} if
$|p-q|>1$ $\forall p\in \alpha$, $q\in \beta$. Suppose $\alpha$,
$\beta_1$, $\ldots$, $\beta_j \subseteq \set{n}$ are sets of
consecutive integers, $\beta_i$ ($i=1,\ldots, j$) are separated, and
\begin{equation} 
 \mbox{ for any } i=1,\ldots, j, \mbox{ there exist } p\in \beta_i
 \mbox{ and } q\in \alpha  \mbox{ such that } |p-q|\leq 1. 
\label{assump}
\end{equation}  
Then $A_{n,k,t}$, $\alpha$, and $\beta\clm=\cup_{i=1}^j \beta_i$
satisfy~Eq.~(\ref{HF}). Indeed,~Eq.~(\ref{HF}) holds for $\alpha$ and
$\beta_1$.  If $1\leq l<j$, then, assuming~Eq.~(\ref{HF}) for $\alpha$ and
$\gamma_l\clm= \cup_{i=1}^l \beta_i$, we have
\begin{eqnarray*}
&&  A_{n,k,t} [\dbl{\alpha}] A_{n,k,t} [\dbl{\gamma_{l+1}}]
=A_{n,k,t} [\dbl{\alpha}] A_{n,k,t} [\dbl{\gamma_l}] A_{n,k,t}
[\dbl{\beta_{l+1}}] \geq \\ && \quad A_{n,k,t} [\dbl{\alpha\cup \gamma_l}]
A_{n,k,t} [\dbl {\alpha \cap \gamma_l}] A_{n,k,t}[\dbl{\beta_{l+1}}].
\end{eqnarray*}
Due to~Eq.~(\ref{assump}), $\alpha \cup \gamma_l$ is a set of consecutive 
integers, so an application of~Eq.~(\ref{HF}) yields
$$ A_{n,k,t} [\dbl{\alpha \cup \gamma_l}] A_{n,k,t}[\dbl{\beta_{l+1}}]
\geq A_{n,k,t} [\dbl{\alpha\cup \gamma_{l+1}}] A_{n,k,t} [\dbl{
(\alpha \cup \gamma_{l})\cap\beta_{l+1}}]. $$
But $(\alpha \cup \gamma_l)\cap \beta_{l+1}=\alpha\cap \beta_{l+1}$
since the sets $\beta_i$ are pairwise disjoint. So, 
\begin{eqnarray} 
 A_{n,k,t}[\dbl{\alpha}] A_{n,k,t}[\dbl{\gamma_{l+1}}] & \geq &
A_{n,k,t}[\dbl{\alpha \cup \gamma_{l+1}}] A_{n,k,t}[\dbl{\alpha
\cap \gamma_l}] A_{n,k,t}[\dbl{\alpha \cap \beta_{l+1}}]= \nonumber \\
&& A_{n,k,t} [\dbl{\alpha \cup \gamma_{l+1}}] A_{n,k,t} [\dbl{\alpha\cap 
\gamma_{l+1}}]. \label{prelim}
\end{eqnarray}
Now, given a set of consecutive integers $\alpha\subseteq \set{n}$ and
any index set $\beta \subseteq \set{n}$, write $\beta=\gamma_1\cup
\gamma_2$ where $\gamma_1\clm =\cup_{i=1}^l\beta_i$,
$\gamma_2\clm=\cup_{i=l+1}^{l+m} \beta_i$, all $\beta_i$ ($i=1,\ldots,
l+m$) are separated, and $\beta_i$ satisfies~Eq.~(\ref{assump}) if and only if $i\leq l$. Then
\begin{eqnarray*}
&& A_{n,k,t} [\dbl{\alpha}] A_{n,k,t}[\dbl{\beta}] =
A_{n,k,t} [\dbl{\alpha}] A_{n,k,t}[\dbl{\gamma_1}] A_{n,k,t}
[\dbl{\gamma_2}] \geq \\
&& \quad A_{n,k,t} [\dbl{\alpha \cup \gamma_1}] A_{n,k,t}[\dbl{\alpha
\cap \gamma_1}] A_{n,k,t}[\dbl{\gamma_2}] = \\
&& \quad A_{n,k,t} [\dbl{\alpha \cup \gamma_1 \cup \gamma_2}] A_{n,k,t} 
[\dbl{\alpha \cap \gamma_1}] = A_{n,k,t} [\dbl{\alpha \cup \beta}] 
A_{n,k,t} [\dbl{\alpha \cap \beta}].
\end{eqnarray*}   
In other words, $A_{n,k,t}$ satisfies~Eq.~(\ref{HF}) if $\alpha\subseteq \set{n}$
is a set of consecutive integers and $\beta \subseteq \set{n}$ is arbitrary.

Finally, if $\alpha_1$, $\alpha_2$, $\beta\subseteq \set{n}$, the sets
$\alpha_i$ ($i=1,2$) are separated, Eq.~(\ref{HF}) holds for $\alpha_1$
and $\beta$, and $\alpha_2$ is a set of consecutive integers, 
then~Eq.~(\ref{HF}) holds for $\alpha\clm =\alpha_1\cup \alpha_2$ and $\beta$:
\begin{eqnarray*}
&& A_{n,k,t}[\dbl{\alpha}] A_{n,k,t}[\dbl{\beta}] =
A_{n,k,t}[\dbl{\alpha_1}] A_{n,k,t}[\dbl{\alpha_2}]
A_{n,k,t}[\dbl{\beta}] \geq \\ && \quad A_{n,k,t} [\dbl{\alpha_1 \cup
\beta}] A_{n,k,t} [\dbl{\alpha_1\cap \beta}] A_{n,k,t}
[\dbl{\alpha_2}] \geq \\ && \quad A_{n,k,t} [\dbl{(\alpha_1 \cup
\beta) \cup \alpha_2}] A_{n,k,t} [\dbl{ (\alpha_1\cup \beta) \cap
\alpha_2}] A_{n,k,t} [\dbl{\alpha_1 \cap \beta}] = \\ && \quad
A_{n,k,t} [\dbl{\alpha \cup \beta}] A_{n,k,t} [\dbl{\alpha_1\cap
\beta}] A_{n,k,t} [\dbl{\alpha_2 \cap \beta}] = \\ && \quad A_{n,k,t}
[\dbl{\alpha \cup \beta}] A_{n,k,t} [\dbl{\alpha \cap \beta}].
\end{eqnarray*}
So, by induction on the number of 'components' of $\alpha$,~Eq.~(\ref{HF}) 
holds for any $\alpha$, $\beta\subseteq \set{n}$. Thus, by the 
Gantmacher-Krein-Carlson theorem, $A_{n,k,t}$ is GKK for any $t\in(0,1)$ 
and any $k$, $n\in \N$.

Now check that $A_{n,k,t}$ have eigenvalue monotonicity if $n\leq
2k+2$ and $t\in(0,1)$ is sufficiently small. Let $\varphi^{k,t}_j
(\lambda)\clm =
\det ( A_{k+j+1,k,t}-\lambda I)$ for $j=1,\ldots, k+1$. Show by 
induction that
\begin{eqnarray}
\varphi^{k,t}_j(\lambda)& = & \cases{ \y^{k+2}-(1-t) & if $j=1$ \cr
\y^{j+k+1}-j(1-t)\y^{j-1}+(j-1)(1-t)^2\y^{j-2}+ & \cr
\frac{t(1-t)^2}{(\y-t)^2}[t^{j-1}-(j-1)t\y^{j-2}+(j-2)\y^{j-1}] & if
$j>1$,}
\label{induct1} \\  
a^{k,t}_j& = &\cases{ (-1)^k(1-t) & if $j=1$ \cr (-1)^{k+j}
t^{j-2}(1-t)^2 & if $j>1$,} \label{induct2} \\
g_j^{k,t}(\lambda)&\clm=& (-1)^{k+1} \det \left(\begin{array}{ccccc} 
\y & 0 & \ldots & 0 & a_j^{k,t} \\ 1 & \y & \ldots & 0 &
a_{j-1}^{k,t} \\ \vdots & \vdots & \ddots & \vdots & \vdots
\\ 0 & 0 & \ldots & \y & a_2^{k,t} \\ 0 & 0 & \ldots & 1 &
a_1^{k,t} \end{array} \right)= \nonumber \\ 
&& -(1-t)\y^{j-1}+(1-t)^2\frac{\y^{j-1}-t^{j-1}}{\y-t} \qquad \qquad 
\qquad \qquad \forall j\in \N.
\label{induct3}
\end{eqnarray}
By direct calculation,
$\varphi_1^{k,t}(\lambda)=\y^{k+2}-(-1)^ka_1^{k,t}$, so, since
$\varphi_1^{k,t}(0)=t$, we have $a_1^{k,t}=(-1)^k(1-t)$.
Thus~Eqs.~(\ref{induct1})--(\ref{induct3}) hold for $j=1$.
Now suppose that $j\geq 2$ and our formulas are true for $j-1$. Expansion of
$\varphi_j(\lambda)$ by its last row gives
\begin{equation}
\varphi_j^{k,t}(\lambda)=\y\varphi_{j-1}^{k,t}(\lambda)+g_j^{k,t}(\lambda).
\label{indphi} \end{equation} 
Since $\varphi_j^{k,t}(0)=t^j$ $\forall j\in \N$, this implies 
$g_j^{k,t}(0)=t^j-t^{j-1}$. On the other hand, expanding $g_j^{k,t}
(\lambda)$ by its first row, we get
\begin{equation}
g_j^{k,t}(\lambda)=\y g_{j-1}^{k,t}(\lambda)+(-1)^{j+k}a_j^{k,t},
\label{indg} \end{equation}
so $a_j^{k,t}=(-1)^{k+j}[g_j^{k,t}(\lambda)-\y 
g_{j-1}^{k,t}(\lambda) ]^{}_{\lambda=0} = (-1)^{j+k}t^{j-2}(1-t)^2$,
which gives~Eq.~(\ref{induct2}). Now, using~Eq.~(\ref{indg}) again together 
with the inductive hypothesis on $g_{j-1}^{k,t}$, we get~Eq.~(\ref{induct3}):
\begin{eqnarray*} 
g_j^{k,t}(\lambda)&=&-(1-t)\y^{j-1}+(1-t)^2\frac{\y^{j-1}-\y t^{j-2}}
{\y-t}+t^{j-2}(1-t)^2= \\
&& -(1-t)\y^{j-1}+(1-t)^2\frac{\y^{j-1}-t^{j-1}}{\y-t}.
\end{eqnarray*}
Finally, substituting the expression for $\varphi_{j-1}^{k,t}(\lambda)$ and
the just verified expression for $g_j^{k,t}(\lambda)$ into~Eq.~(\ref{indphi})
yields~Eq.~(\ref{induct1}).  

If $\set{n}\supseteq \alpha=\cup_{i=1}^j \alpha_i$ is the union of 
separated sets of consecutive integers, then 
$ \det (A_{n,k,t}(\dbl{\alpha})
-\lambda I)= \prod_{i=1}^j \det (A(\dbl{\alpha_i})-\lambda I) $ 
since $A_{n,k,t}-\lambda I$ is Hessenberg (the same observation earlier 
led to~Eq.~(\ref{*})). Since $A_{n,k,t}-\lambda I$ is Toeplitz, the  product
in the right hand side equals $\prod_{i=1}^j \det (A_{n,k,t} (\dbl
{\set{\# \alpha_i}}) -\lambda I)$. Hence, to prove eigenvalue monotonicity 
of $A_{n,k,t}$ for $n\leq 2k+2$ it is enough to prove it for
leading principal submatrices of $A_{n,k,t}$ only, i.e., to show
$$ l(A_{k+j+1,k,t})\leq l(A_{k+j,k,t}) \qquad \forall j\in \N,$$ i.e.,
that $\varphi_j^{k,t}$ has a root in $(0,1]$ for any $j\leq k+1$, and
$$\min\{ \lambda\in(0,1] : \varphi_j(\lambda)=0\} \leq \min\{ \lambda
\in (0,1] : \varphi_{j-1}(\lambda)=0\}, \qquad j=2,\ldots, k+1 $$
(since $A_{k+j,k,t}$ is a $P$-matrix, the coefficients of its characteristic
polynomial are strictly alternating, so $A_{k+j,k,t}$ has no nonpositive
eigenvalues).  Observe that
$\varphi_j^{k,t}(\lambda)=t^j-\lambda\widetilde{\varphi}_j^{k,t}
(\lambda)$ where
\begin{eqnarray*} 
\widetilde{\varphi}_j^{k,t}(0)& = & -\frac{d\varphi_j^{k,t}(\lambda)}
{d\lambda}\bigg|_{\lambda=0}
\belowrightarrow{t\to 0+}=-\frac{d\nu_j^k(\lambda)}{d\lambda}\bigg|_{\lambda=0}, \\
\nu_j^k(\lambda)&\clm=&\lim_{t\to 0+} \varphi_j^{k,t}(\lambda)=
\y^{j+k+1}-j\y^{j-1}+(j-1)\y^{j-2}. 
\end{eqnarray*}   
So, $\lim_{t\to 0+} \widetilde{\varphi}_j^{k,t}(0)=k+3-j\geq 2$
$\; \forall j=1,\ldots, k+1$.

Since $0$, the minimal real root of $\nu_j^k$, is simple, the
minimal real root $\lambda_j$ of $\varphi_j^{k,t}$ is positive and
simple for all $j=1,\ldots, k+1$ whenever $t$ is sufficiently
small. But then $\tilde{\varphi}_j^{k,t}(\lambda_j)$ is bounded below
by a positive constant for any $j=1,\ldots, k+1$, hence
$$ \lambda_j=\frac{t^j}{\widetilde{\varphi}_j^{k,t}(\lambda_j)} <
\frac{t^{j-1}}{\widetilde{\varphi}_{j-1}^{k,t}(\lambda_{j-1})}=\lambda_{j-1}
\qquad \forall j=1,\ldots, k+1$$ if $t$ is small. So, for any $k\in \N$ 
and $n\leq 2k+2$, there exists $t(k)\in (0,1)$ such that $A_{n,k,t}$ 
is a $\tau$-matrix for all $t\in (0,t(k))$.

Now let $B_k\clm=\lim_{t\to 0+} A_{2k+2,k,t}$. The matrix $B_k$ is
Toeplitz with first column $$(1,1,\underbrace{0, \ldots, 0}_{2k \; {\rm
\scriptstyle times}})^T$$ and first row $$(1,\underbrace{ 0, \ldots,
0}_{k\; {\rm \scriptstyle times}}, (-1)^k, (-1)^k,
\underbrace{0,\ldots, 0}_{k-1 \; {\rm \scriptstyle times}}).$$
Show that there exists $K\in \N$ such that, for all $k>K$, $B_k$ has
an eigenvalue $\lambda$ with $\Re \lambda<0$.  As the eigenvalues
depend continuously on the entries of the matrix, this will
demonstrate that, for any $k>K$, there exists $t\in (0,1)$
such that the GKK $\tau$-matrix $A_{2k+2,k,t}$ has an
eigenvalue with negative real part.

The polynomial $\nu_{k+1}^k$ has a root with negative real part iff
the polynomial $\psi_k$ where
$$\psi_k(\lambda)\clm=\frac{\nu^k_{k+1}(-\lambda)}{(1+\lambda)^{k-1}}=
(1+\lambda)^{k+3}-(k+1)(1+\lambda)+k$$ has a root with positive real
part. Since
$$\psi_k(\lambda)=\lambda \biggl[ \sum_{j=0}^{k+1} {k+3 \choose j}
\lambda^{k+3-j-1}+2 \biggr],$$ it is, in turn, enough to show that
$\eta_k$ where
$$\eta_k(\lambda)\clm =\lambda^{k+3}\psi_k\left(\frac{1}{\lambda}\right)=
2\lambda^{k+2}+\sum_{j=2}^{k+3} {k+3 \choose j} \lambda^{k+3-j}$$ has
a root with positive real part. The Hurwitz matrix for the polynomial
$\eta_k$ is
$$ H_k\clm= \left( \begin{array}{cccccc} {k+3 \choose 2} & {k+3
\choose 4} & {k+3 \choose 6} & {k+3 \choose 8} & {k+3 \choose 10} &
\cdots \mystrut \\ 2 & {k+3 \choose 3} & {k+3 \choose 5} & {k+3
\choose 7} & {k+3 \choose 9} & \cdots \mystrut \\ 0 & {k+3 \choose 2}
& {k+3 \choose 4} & {k+3 \choose 6} & {k+3 \choose 8} & \cdots
\mystrut \\ 0 & 2 & {k+3 \choose 3} & {k+3 \choose 5} & {k+3 \choose
7} & \cdots \mystrut \\ 0 & 0 & {k+3 \choose 2} & {k+3 \choose 4} &
{k+3 \choose 6} & \cdots \mystrut \\ \vdots & \vdots & \vdots & \vdots
& \vdots & \ddots \mystrut
\end{array} \right)_{(k+2)\times (k+2)} \goback{$\scriptstyle (k+2)(k+2)$}. $$
Compute the minor $H_k[\dbl{2\clm5}]$, taking out the factors $k+3 \choose 2$, 
$k+3 \choose 4$, $k+3 \choose 6$ from its second, third, and fourth columns
respectively. We obtain
$$ H_k[\dbl{2\clm5}]=-\frac{1}{132300}(3k^3-49k^2-210k-318)(k+4)^2(k+5)
{k+3 \choose 2}{k+3 \choose 4} {k+3\choose 6}.$$ It follows that
$H_k[\dbl{2\clm5}]<0$ for $k$ large enough, precisely, for all
$k>20$. But the Hurwitz matrix of a nonpositive stable polynomial is
totally nonnegative (see~\cite{4}). So, for $k>20$, $\eta_k$ has a
zero with positive real part, therefore, $\nu^k_{k+1}$ has a zero 
with negative real part. This completes the proof that the GKK $\tau$-matrices
$A_{2k+2,k,t}$ are unstable for sufficiently large $k$ and small $t$.

\paragraph{Remark} {\em {
To illustrate the result, consider the matrix $A_{44,21,1/2}$, i.e., 
the Toeplitz matrix whose first column is $$(1,1, \underbrace{0,\ldots,
0}_{42 \; {\rm \scriptstyle times}})^T$$ 
and first row is
$$ (1, \underbrace{0,\ldots, 0}_{21 \; {\rm \scriptstyle times}}, -1/2,
-1/2^2, 1/2^3, -1/2^4, \ldots, -1/2^{22})$$
 and the limit matrix
$B_{21}$, with the same first column as $A_{44,21,1/2}$ and first 
row equal to $$(1, \underbrace{0,\ldots,0}_{21 \; 
{\rm \scriptstyle times}}, -1, -1, \underbrace{0, \ldots, 
0}_{20 \; {\rm \scriptstyle times}}).$$
According to MATLAB,
the two eigenvalues with minimal real part of the first matrix are
$$-2.809929189497896\cdot 10^{-2} \pm 3.275076252367531\cdot 10^{-1}i;$$ 
those of the second are $$-3.420708309454068\cdot 10^{-2} \pm
3.400425852703498\cdot 10^{-1}i.$$}}
\hspace{0.6cm} 

\section*{Acknowledgements}

 I am grateful to Professor Hans Schneider for many valuable remarks
and suggestions and to the referee for careful reading.

\end{document}